\documentclass[a4paper, 12pt]{amsart}
\usepackage{amssymb, latexsym, amsmath}

\theoremstyle{plain}
\newtheorem{theorem}{Theorem}[section]
\newtheorem{lemma}[theorem]{Lemma}
\newtheorem{corollary}[theorem]{Corollary} 

\theoremstyle{definition}

\theoremstyle{remark}
\newtheorem*{remark}{Remark}

\numberwithin{equation}{section} 

\begin {document}
\baselineskip=24pt
\title[Finite complete rewriting systems]{On finite complete rewriting systems and large subsemigroups}

\author {K.B. Wong}
\address {Institute of Mathematical Sciences \\
University of Malaya \\
50603 Kuala Lumpur, Malaysia}
\email {kbwong@um.edu.my}

\author {P.C. Wong}
\address {Institute of Mathematical Sciences \\
University of Malaya \\
50603 Kuala Lumpur, Malaysia}
\email {wongpc@um.edu.my}


\begin {abstract} Let $S$ be a semigroup and $T$ be a subsemigroup of finite index in $S$ (that is, the set $S\setminus T$ is finite). The subsemigroup $T$ is also called a large subsemigroup of $S$. It is well known that if $T$ has a finite complete rewriting system, then so does $S$. In this paper, we will prove the converse, that is, if $S$ has a finite complete rewriting system, then so does $T$. Our proof is purely combinatorial and also constructive.
\end {abstract}

\maketitle

\section {\bf Introduction.}

Let $S$ be a semigroup and $T$ be a subsemigroup of finite index in $S$ (that is, the set $S\setminus T$ is finite). Then $T$ is called a \emph {large subsemigroup} of $S$, and $S$ is called a \emph {small extension} of $T$. In \cite{Ruskuc}, Ru\v{s}kuc asked if $S$ is a small extension of $T$, whether $S$ has a finite complete rewriting system if and only if $T$ has a finite complete rewriting system (see \cite[Problem 11.1 (iii)]{Ruskuc} and \cite[Remark 4.2]{Wang2}). This problem was partially solved by Wang in \cite[Theorem 1]{Wang}, who proved that if $T$ has a finite complete rewriting system, then so does $S$. However it is still not known whether $T$ has a finite complete rewriting system or not, when $S$ has a finite complete rewriting system. In this paper we shall prove that this is true, i.e., we shall prove the following:

\begin{theorem}\label{main_theorem} Suppose $S$ is a small extension of $T$. If $S$ has a finite complete rewriting system, then so does $T$.
\end{theorem}

By Theorem \ref{main_theorem} and the result of Wang \cite[Theorem 1]{Wang}, we have completely answered the problem posed by Ru\v{s}kuc (see \cite[Problem 11.1 (iii)]{Ruskuc}).

\begin{corollary}\label{consequence_main_theorem} Suppose $S$ is a small extension of $T$. Then $S$ has a finite complete rewriting system if and only if $T$ has a finite complete rewriting system.
\end{corollary}

Let $A$ be a non-empty set. This set $A$ is called the alphabet and the elements of $A$ are called letters. We shall denote the free semigroup and free monoid on $A$ by $A^+$ and $A^*$, respectively. The elements of $A^+$ and $A^*$ are called words. Note that $A^*=A^+\cup \{\epsilon\}$, where $\epsilon$ is the empty word. Given a word  $W\in A^*$, we shall denote its length by $\Vert W\Vert$, defined as the numbers of letters in $W$.

A \emph{rewriting system} $R$ over $A$ is a set of rules $U\to V$, which are elements of $A^+\times A^+$. A word $W_1\in A^+$ is said to be rewritten to another word $W_2\in A^+$ by a \emph{one-step reduction} induced by $R$, if $W_1=Z_1XZ_2$ and $W_2=Z_1YZ_2$ for some rule $X\to Y$ in $R$. In this situation we write $W_1\to_R W_2$. The reflexive transitive closure and the reflexive symmetric transitive closure of $\to_R$ are denoted by $\to_R^*$ and $\leftrightarrow_R^*$, respectively. The relation $\leftrightarrow_R^*$ is defined to be the congruence on $A^+$ generated by $R$ and it defines the quotient semigroup $S=A^+/\leftrightarrow_R^*$. $S$ is said to be presented by the \emph{semigroup presentation} $[A ; R]$. If both $A$ and $R$ are finite, we say the semigroup presentation is finitely presented. For $U\in A^+$, $[U]_R$ shall denote the class of $U$ modulo $\leftrightarrow_R^*$.

Let $\textnormal {Left} (R)=\{ X\in A^+\ :\ X\to Y\in R\}$ and $\textnormal {Irr} (R)=A^+\setminus A^*\textnormal {Left} (R)A^*$. Obviously,  $\textnormal {Irr} (R)$ is the set of all words in $A^+$ that cannot be reduced by  any rule in $R$. A word $W\in A^+$ is called an \emph{irreducible word} if $W\in \textnormal {Irr} (R)$.

We say $R$ is \emph{Noetherian} if there is no infinite reduction sequence,
\begin{equation}
W_1\to_R W_2\to_R W_3\to_R\cdots .\notag
\end{equation}
$R$ is said to be \emph{confluent} if whenever $U\to_R^* V$ and $U\to_R^* W$, then there is an $X\in A^+$ such that $V\to_R^* X$ and $W\to_R^* X$. If $R$ is both Noetherian and confluent, we say that $R$ is a \emph{complete rewriting system}. 

The following fact is well known.

\begin{theorem}\label{complete_rewriting_Irr} Suppose $R$ is a complete rewriting system. Then for each $W\in A^+$, there is a unique $W'\in \textnormal {Irr} (R)$ such that $W\to_R^* W'$.
\end{theorem}

Theorem \ref{complete_rewriting_Irr} will be used implicitly in many parts of the paper.
 Let $R$ be a complete rewriting system on $A^+$. Then given any word $W\in A^+$, by Theorem \ref{complete_rewriting_Irr}, there is a $U\in \textnormal{Irr}(R)$ such that 
 
\begin{equation}
W\to_RW_1\to_RW_2\to_R\cdots\to_RW_n=U.\notag
\end{equation}

The length of the above reduction sequence starting with $W$ and ends with $U$ is $n$. The \emph{disorder} of $W$, denoted by $d_R(W)$, is the maximum of the lengths of all of the reduction sequences starting with $W$ and ends with $U$. Note that $d_R(W)$ is finite. Suppose it is not. Then there is a $V_1\in A^+$ such that $W\to_RV_1$ and $d_R(V_1)$ is infinite, for the number of subwords of $W$ that are contained in $\textnormal {Left} (R)$ is finite. Then again there is a $V_2\in A^+$ such that $V_1\to_RV_2$ and $d_R(V_2)$ is infinite, and this process can go on indefinitely. So $W\to_R V_1\to_R V_2\to_R\cdots$ is an infinite reduction sequence, a contradiction. Note also that $W\in \textnormal{Irr}(R)$ if and only if $d_R(W)=0$ (see \cite{GS} and \cite{HM}).

The following useful lemma is obvious.

\begin{lemma}\label{disorder_inequality} If $U\to_R V$, then $d_R(U)>d_R(V)$. Furthermore if $W$ is a subword of $U$, then $d_R(U)\geq d_R(W)$.
\end{lemma}

A semigroup is said to have a \emph{finite complete rewriting system} if it has a finitely presented semigroup presentation for which the rewriting system is complete.

\section {\bf A criterion.}

Let $[ A\ ;\ R]$ be a finitely presented semigroup presentation for $S$ for which $R$ is complete. Let $T$ be a subsemigroup of $S$. In this section we first prove a criterion for $[B\ ;\ R_T]$ to be a semigroup presentation for $T$ where $B$ is any non-empty set and $R_T$ is a complete rewriting system over $B$. This will be done in Theorem \ref{general_propert_R}. Then by replacing $T$ with $S$ we can get $[B\ ;\ R_S]$ to be a semigroup presentation for $S$ and $R_S$ is a complete rewriting system over $B$. This will be done in Corollary \ref{general_propert_presentation_S}

Let $A(T)$ be a subset of $A^+$ such that
\begin{equation}
\{ W\in \textnormal{Irr} (R)\ :\ [W]_R\in T\}\subseteq A(T)\subseteq \{ W\in A^+\ :\ [W]_R\in T\}.\notag
\end{equation}
Let $(B,R_T,A(T),\phi, \rho)$ be a  5-tuple where $B$ is a non-empty set, $R_T$ is a rewriting system over $B$, $\phi : B^+\to A^+$ is a homomorphism with $[\phi(U')]_R\in T$ for all $U'\in B^+$, and $\rho : A(T)\to B^+$ is a function. We say the 5-tuple $(B,R_T,A(T),\phi, \rho)$  has \emph{Property} $\mathcal R$ relative to $[ A\ ;\ R]$, if it satisfies the following:
\begin{itemize}
\item [(P1)] for any $U\in A(T)$ and $V_1\in A^+$ with $U\to_R V_1$, there is a $U'\in B^+$ such that $U\to_R V_1\to_R^* \phi(U')$ and  $\rho(U)\to_{R_T} U'$,
\item [(P2)] for any $U',V'\in B^+$ with $U'\to_{R_T}^* V'$, we have $\phi(U')\to_{R}^* \phi(V')$,
\item [(P3)] there does not exist an infinite reduction sequence 
\begin{equation}
U_1'\to_{R_T}U_2'\to_{R_T}U_3'\to_{R_T}\cdots,\notag
\end{equation}
 of words from $B^+$ such that $\phi(U_1')=\phi(U_2')=\phi(U_3')=\cdots$,
\item [(P4)] for each $U'\in B^+$ there is a $U''\in B^+$ such that $\phi(U'')\in A(T)$ and $U'\to_{R_T}^* U''$,
\item [(P5)] $\phi(\rho(U))=U$ for all $U\in A(T)$,
\item [(P6)] $U'\to_{R_T}^*\rho(\phi (U'))$ for all $U'\in B^+$ with $\phi (U')\in A(T)$.
\end{itemize}

\begin{lemma}\label{condition_P11} Suppose \textnormal{(P1)}, \textnormal{(P2)}, \textnormal{(P4)}, and \textnormal{(P6)} hold. Then for any $U\in A(T)$ and $V\in\textnormal {Irr} (R)$ with $U\to_R^* V$, we have $V\in A(T)$ and $\rho(U)\to_{R_T}^* \rho(V)$.
\end{lemma}

\begin{proof} By the definition of $A(T)$, clearly $V\in A(T)$. We shall prove by induction on $d_R(U)$ that $\rho(U)\to_{R_T}^* \rho(V)$.

Suppose $d_R(U)=0$ then $U=V$. Thus $\rho(U)=\rho(V)$ and $\rho(U)\to_{R_T}^* \rho(V)$. Suppose $d_R(U)>0$. Assume that it is true for all $U_1$ with $d_R(U_1)<d_R(U)$.

Let $U\to_R V_1\to_R^* V$. By (P1), there is a $U'\in B^+$ such that $U\to_R V_1\to_R^* \phi(U')$ and  $\rho(U)\to_{R_T} U'$. By (P4), there is a $U''\in B^+$ such that $\phi(U'')\in A(T)$ and $U'\to_{R_T}^* U''$. By (P2), $\phi(U')\to_{R}^* \phi(U'')$. Therefore $U\to_R^* \phi(U'')$ and $\rho(U)\to_{R_T}^* U''$. Since $V\in\textnormal{Irr} (R)$, we have $\phi(U'')\to_R^* V$. Furthermore $d_R(\phi(U''))<d_R(U)$ (by Lemma \ref{disorder_inequality}). Therefore by induction $\rho(\phi(U''))\to_{R_T}^* \rho(V)$. Now by (P6), $U''\to_{R_T}^* \rho(\phi(U''))$. Hence $\rho(U)\to_{R_T}^* \rho(V)$.

The proof of this lemma is complete.
\end{proof}

\begin{theorem}\label{general_propert_R} If $(B,R_T,A(T),\phi,  \rho)$ has \emph{Property} $\mathcal R$ relative to $[ A\ ;\ R]$, then $[B\ ;\ R_T]$ is a semigroup presentation for $T$ and $R_T$ is complete.
\end{theorem}

\begin{proof} We will first prove that $[B\ ;\ R_T]$ is a semigroup presentation for $T$. Let $\psi : [B\ ;\ R_T]\to T$ be defined by $\psi([U']_{R_T})=[\phi(U')]_R$ for all $U'\in B^+$. Now we show that $\psi$ is well defined. It is sufficient to prove $U'\to_{R_T} V'$ for $V'\in B^+$ implies that $\phi(U')\to_R^* \phi(V')$. This fact follows from (P2), so $\psi$ is well defined.

Now we show that $\psi$ is a homomorphism. Let $U',V'\in B^+$. Then $\psi([U'V']_{R_T})=[\phi(U'V')]_R=[\phi(U')\phi(V')]_R=[\phi(U')]_R[\phi(V')]_R=\psi([U']_{R_T})\psi([V']_{R_T})$, where the second equality follows from the fact that $\phi$ is a homomorphism.

Now we show that $\psi$ is surjective. Let $[W]_R\in T$ for some $W\in A^+$. Since $R$ is complete, we may assume $W\in\textnormal {Irr} (R)$. Note that $W\in A(T)$, so $\psi([\rho(W)]_{R_T})=[\phi(\rho(W))]_R=[W]_R$, where the last equality follows from (P5). Hence $\psi$ is surjective. 

Now we show that $\psi$ is injective. Let $U',V'\in B^+$ with $\psi([U']_{R_T})=\psi([V']_{R_T})$. Then $[\phi(U')]_R=[\phi(V')]_R$. By (P4), there are $U'', V''\in B^+$ such that $\phi (U''), \phi(V'')\in A(T)$, $U'\to_{R_T}^* U''$, $V'\to_{R_T}^* V''$. By (P2),  $\phi(U')\to_{R}^* \phi(U'')$ and $\phi(V')\to_{R}^* \phi(V'')$. So $[\phi(U'')]_R=[\phi(V'')]_R$. Since $R$ is complete, there is a $V_1\in \textnormal {Irr} (R)$ such that $\phi(U'')\to_R^* V_1$ and $\phi(V'')\to_R^* V_1$. By Lemma \ref{condition_P11}, $\rho(\phi (U''))\to_{R_T}^* \rho (V_1)$, and then by (P6), $U''\to_{R_T}^* \rho (V_1)$. Therefore $U'\to_{R_T}^* \rho (V_1)$. Similarly, we have $V'\to_{R_T}^* \rho (V_1)$. Hence  $[U']_{R_T}=[\rho(V_1)]_{R_T}=[V']_{R_T}$ and $\psi$ is injective.

Now we have shown that $[B\ ;\ R_T]$ is a semigroup presentation for $T$, via $\psi$. We will now proceed to prove that $R_T$ is  complete.  

Suppose $R_T$ is not Noetherian. Then there exists an infinite reduction sequence 
\begin{equation}
U_1'\to_{R_T}U_2'\to_{R_T}U_3'\to_{R_T}\cdots,\notag
\end{equation}
 of words from $B^+$. By (P2), $\phi(U'_i)\to_{R}^* \phi(U'_{i+1})$ for all $i$. Since $R$ is Noetherian, there is an integer $i_0$ such that for all $i\geq i_0$,  $\phi(U'_i)=\phi(U'_{i+1})$, but this contradicts (P3). Hence $R_T$ is Noetherian.

Now we prove that $R_T$ is confluent. Suppose $U'\to_{R_T}^* V_1'$ and $U'\to_{R_T}^* V_2'$ with $U',V_1',V_2'\in B^+$. By (P4), we may assume $\phi (V_1'), \phi (V_2')\in A(T)$. Since $R$ is complete, there is a $V_3\in \textnormal {Irr} (R)$ with $\phi (V_1')\to_R^* V_3$ and $\phi(V_2')\to_R^* V_3$. By Lemma \ref{condition_P11}, $\rho(\phi(V_1'))\to_{R_T}^* \rho (V_3)$, and then by (P6) $V_1'\to_{R_T}^* \rho (V_3)$. Similarly $V_2'\to_{R_T}^* \rho (V_3)$. Hence $R_T$ is confluent and is complete.
\end{proof}

In the case when $T=S$ and there is a 5-tuple $(B,R_S,A(S),\phi, \rho)$  that has \emph{Property} $\mathcal R$ relative to $[ A\ ;\ R]$, we have the following corollary:

\begin{corollary}\label{general_propert_presentation_S} $[B\ ;\ R_S]$ is a semigroup presentation for $S$ and $R_S$ is  complete.
\end{corollary}

\section {\bf Changing the semigroup presentation for $S$.}

Let $[ A\ ;\ R]$ be a finitely presented semigroup presentation for $S$ for which $R$ is complete. Let $W_0\in A^+$ be such that $\Vert W_0\Vert>1$ and $W_0\in \textnormal {Irr} (R)$. Now let $s$ be a letter that does not appear in $A$ and set $B=A\cup \{s\}$. We wish to find a complete rewriting system $R_S$ such that $[B\ ;\ R_S]$ is a finitely presented semigroup presentation for $S$ and $W_0\to_{R_S}^* s$.

By Corollary \ref{general_propert_presentation_S}, we need to find a 5-tuple $(B,R_S,A(S),\phi,  \rho)$  that has Property $\mathcal R$ relative to $[ A\ ;\ R]$. Note that $B$ has been defined and is finite.
 
Let $A(S)=A^+$. Let $\phi_1: B\to A^+$ be defined by $\phi_1(a)=a$ for all $a\in A$ and $\phi_1(s)=W_0$. Clearly $\phi_1$ can be extended to a homomorphism $\phi : B^+\to A^+$ by defining $\phi (U')=\phi_1(b_1)\dots \phi_1(b_l)$ for all $U'=b_1\dots b_l\in B^+$. For convenience, we may define $\phi(\epsilon_B)=\epsilon_A$ where $\epsilon_B$ and $\epsilon_A$ are empty words in $B^*$ and $A^*$, respectively.

Recall that we have set $A(S)=A^+$. We define $\rho : A(S)\to B^+$ as follows :

Let $W\in A(S)$. 
\begin{itemize}
\item [(a)] If $W$ ends with the subword $W_0$, say $W=X_1W_0$ for some $X_1\in A^*$ (we use $A^*$ instead of $A^+$ because we allow $X_1$ to be the empty word), then $\rho (W)=\rho (X_1)s$ (in the event $X_1=\epsilon_A$,  set $\rho (W)=s$). 
\item [(b)] Suppose $W$ does not end with the subword $W_0$. Let $W=X_2a$ for some $X_2\in A^*$ and $a\in A$. Set $\rho (W)=\rho (X_2)a$ (in the event $X_2=\epsilon_A$, set $\rho (W)=a$). 
\end{itemize}

As for the homomorphism $\phi$, we may define $\rho(\epsilon_A)=\epsilon_B$.

\begin{lemma}\label{S_property_rho_ending} Let $X_1, X_2, X_3\in A(S)$. If $\rho (X_1X_2X_3)=\rho(X_1X_2)\rho(X_3)$, then $\rho (X_2X_3)=\rho(X_2)\rho(X_3)$.
\end{lemma}
\begin{proof} We prove by induction on $\Vert X_3\Vert$. Clearly it holds if $\Vert X_3\Vert=0$, i.e., $X_3$ is the empty word. Suppose $\Vert X_3\Vert>0$. Assume that it holds for all $X_4$ with $\Vert X_4\Vert <\Vert X_3\Vert$.

\noindent
{\bf Case 1.} Suppose $X_3$ ends with the subword $W_0$, say $X_3=X_4W_0$ for some $X_4\in A^*$. Then $\rho(X_1X_2X_3)=\rho(X_1X_2X_4)s$, $\rho(X_2X_3)=\rho(X_2X_4)s$ and $\rho(X_3)=\rho(X_4)s$. Since $\rho (X_1X_2X_3)=\rho(X_1X_2)\rho(X_3)$, we have $\rho(X_1X_2X_4)=\rho(X_1X_2)\rho(X_4)$. By induction, $\rho (X_2X_4)=\rho(X_2)\rho(X_4)$. Therefore $\rho(X_2X_3)=\rho(X_2X_4)s=\rho(X_2)\rho(X_4)s=\rho(X_2)\rho(X_3)$.

 \noindent
{\bf Case 2.} Suppose $X_3$ does not  end with the subword $W_0$. Let $X_3=X_4a$ for some $a\in A$ and $X_4\in A^*$. Then $\rho(X_3)=\rho(X_4)a$. Now $\rho (X_1X_2X_3)=\rho(X_1X_2)\rho(X_3)=\rho(X_1X_2)\rho(X_4)a$. So $\rho (X_1X_2X_3)$ is a word in $B^+$ that ends with the letter $a$.

We claim that $X_1X_2X_3$ does not end with the subword $W_0$. Suppose the contrary. Then $X_1X_2X_3=Z_1W_0$ for some $Z_1\in A^*$ and $\rho (X_1X_2X_3)=\rho(Z_1)s$. So $\rho (X_1X_2X_3)$ is a word in $B^+$ that ends with the letter $s$. But this contradicts the last sentence of the previous paragraph. Thus our claim has been established.
Therefore $X_1X_2X_3=X_1X_2X_4a$ and $\rho(X_1X_2X_3)=\rho(X_1X_2X_4)a$. This implies that $\rho(X_1X_2X_4)=\rho(X_1X_2)\rho(X_4)$, and by induction $\rho(X_2X_4)=\rho(X_2)\rho(X_4)$.

Note also that $X_2X_3$ does not end with the subword $W_0$, for otherwise $X_1X_2X_3$ would end with the subword $W_0$. Therefore $X_2X_3=X_2X_4a$ and $\rho(X_2X_3)=\rho(X_2X_4)a$. Since $\rho(X_2X_4)=\rho(X_2)\rho(X_4)$ and $\rho(X_3)=\rho(X_4)a$, we conclude that $\rho (X_2X_3)=\rho(X_2)\rho(X_3)$.
\end{proof}

\begin{lemma}\label{S_property_rho} Let $X_1, X_2\in A(S)$. Then  either
\begin{itemize}
\item [(a)]  $\rho(X_1X_2)=\rho(X_1)\rho(X_2)$, or 
\item [(b)] $\rho(X_1X_2)=\rho(Z_1)s\rho(Z_4)$ where $X_1=Z_1Z_2$, $X_2=Z_3Z_4$ and $Z_2Z_3=W_0$ ($Z_1,Z_4\in A^*$ and $Z_2,Z_3\in A^+$).
\end{itemize}
\end{lemma}
\begin{proof} We prove by induction on $\Vert X_2\Vert$.  Clearly it holds if $\Vert X_2\Vert=0$, i.e., $X_2$ is the empty word. Suppose $\Vert X_2\Vert>0$. Assume that it holds for all $X_3$ with $\Vert X_3\Vert <\Vert X_2\Vert$.

\noindent
{\bf Case 1.} Suppose $X_2$ ends with the subword $W_0$, say $X_2=X_3W_0$ for some $X_3\in A^*$. Then $\rho(X_1X_2)=\rho(X_1X_3)s$. If $X_3$ is the empty word, then $\rho(X_1X_2)=\rho(X_1)s=\rho(X_1)\rho(X_2)$, we are done. If $X_3$ is not the empty word, then $\rho(X_1X_2)=\rho(X_1X_3)s$, and by induction ($\Vert X_3\Vert<\Vert X_2\Vert$), either
 $\rho(X_1X_3)=\rho(X_1)\rho(X_3)$ or $\rho(X_1X_3)=\rho(Z_1)s\rho(Z_4)$, where  $X_1=Z_1Z_2$, $X_3=Z_3Z_4$ and $Z_2Z_3=W_0$ ($Z_1,Z_4\in A^*$ and $Z_2,Z_3\in A^+$). Suppose the former holds. Then $\rho(X_2)=\rho(X_3W_0)=\rho(X_3)s$ and $\rho(X_1X_2)=\rho(X_1X_3)s=\rho(X_1)\rho(X_3)s=\rho(X_1)\rho(X_2)$.
 
 Suppose the latter holds. Then $X_2=Z_3Z_4W_0=Z_3Z_5$ ($Z_5=Z_4W_0$) and $\rho(X_1X_2)=\rho(X_1X_3)s=\rho(Z_1)s\rho(Z_4)s=\rho(Z_1)s\rho(Z_4W_0)=\rho(Z_1)s\rho(Z_5)$. Thus  the lemma holds.
 
 \noindent
{\bf Case 2.} Suppose $X_2$ does not  end with the subword $W_0$ but $X_1X_2$ ends with the subword $W_0$, say $X_1X_2=X_3W_0$ for some $X_3\in A^*$.  Then $\Vert W_0\Vert >\Vert X_2\Vert$ and $X_1=X_3X_4$ where $X_4X_2=W_0$ ($X_4\in A^+$). Note that $\rho(X_1X_2)=\rho(X_3)s$ and the lemma holds.

\noindent
{\bf Case 3.} Suppose $X_1X_2$ does not  end with the subword $W_0$. Let $X_2=X_3a$ where $a\in A$ and $X_3\in A^*$. Then $\rho(X_1X_2)=\rho(X_1X_3)a$. Since $\Vert X_3\Vert<\Vert X_2\Vert$, by induction and using an argument similar to Case 1, we conclude that the lemma holds.
\end{proof}

\begin{lemma}\label{S_condition_P5} $\phi(\rho(U))=U$ for all $U\in A(S)$. \textnormal{(Property (P5)).}
\end{lemma}
\begin{proof} Let $U\in A(S)$. We shall prove by induction on $\Vert U\Vert$ that $\phi(\rho(U))=U$. If $\Vert U\Vert=1$, then $U=a$ for some $a\in A$ and clearly $\phi(\rho(U))=a=U$. Suppose $\Vert U\Vert>1$. Assume the lemma holds for all $U_1\in A(S)$ with $\Vert U_1\Vert<\Vert U\Vert$. 

Suppose $U$ ends with the subword $W_0$, say $U=X_1W_0$ for some $X_1\in A^*$. Then $\rho (U)=\rho (X_1)s$ and $\phi(\rho (U))=\phi(\rho (X_1))\phi(s)=\phi(\rho (X_1))W_0=X_1W_0=U$, where the first equality follows from the fact that $\phi$ is a homomorphism, and the second last equality follows from induction (clearly $\Vert X_1\Vert<\Vert U\Vert$).

Suppose $U$ does not end with the subword $W_0$. Let $U=X_2a$ for some $X_2\in A^*$ and $a\in A$. Now $\rho (U)=\rho (X_2)a$ and  similarly by induction $\phi(\rho (U))=\phi(\rho (X_2))\phi(a)=\phi(\rho (X_2))a=X_2a=U$. Hence the lemma holds.
\end{proof}

Now we define the rules in $R_S$. Recall that $W_0\in \textnormal {Irr} (R)$ and $\Vert W_0\Vert>1$.

\begin{itemize}
\item [($\mathcal C$1)] for each $X\to Y\in R$ put $\rho (X)\to \rho (Y)$ in $R_S$;
\item [($\mathcal C$2)] put $W_0\to s$ in $R_S$;
\item [($\mathcal C$3)] if there is a rule $X_1X_2\to Y_1\in R$ such that $W_0=Z_1X_1$ ($X_1,Y_1\in A^+$ and $X_2,Z_1\in A^*$), put $\rho (Z_1X_1X_2)\to \rho (Y')$ in $R_S$ where $Z_1X_1X_2\to_R^* Y'$ and $Y'\in\textnormal {Irr} (R)$;
\item [($\mathcal C$4)] if there is a rule $X_2X_1\to Y_1\in R$ such that $W_0=X_1Z_1$ ($X_1,Y_1\in A^+$ and $X_2,Z_1\in A^*$), put $\rho (X_2X_1Z_1)\to \rho (Y')$ in $R_S$ where $X_2X_1Z_1\to_R^* Y'$ and $Y'\in\textnormal {Irr} (R)$;
\item [($\mathcal C$5)] if there is a rule $X_2X_3X_4\to Y_1\in R$ such that $W_0=X_4X_5=X_1X_2$ ($X_2, X_4, Y_1\in A^+$ and $X_3,X_5,X_1\in A^*$), put $\rho (X_1(X_2X_3X_4)X_5)\to \rho (Y')$ in $R_S$ where 
 $X_1(X_2X_3X_4)X_5\to_R^* Y'$ and $Y'\in\textnormal {Irr} (R)$;
\item [($\mathcal C$6)] if there are $X_1, X_2, X_3\in A^+$ such that $W_0=X_1X_2=X_2X_3$, put $sX_3\to X_1s$ in $R_S$ (in the event of this we must have $\Vert X_1\Vert=\Vert X_3\Vert$).
\end{itemize}
Note that the number of rules  of the form $\mathcal C$1 and $\mathcal C$2 that we put in $R_S$ is finite. The number of rules of the form $\mathcal C$3 that we put in $R_S$ is also finite because $R$ is finite and $W_0$ is a fixed word. Similarly for the number of rules of the form $\mathcal C$4 up to  $\mathcal C$6. Therefore $R_S$ is a finite rewriting system. 

\begin{remark} Note that one can subsume the rules ($\mathcal C$1), ($\mathcal C$3), and ($\mathcal C$4) all within ($\mathcal C$5) by just allowing $X_1X_2$ and $X_4X_5$ be empty, as well as equal to $W_0$.
\end{remark}

Since $A(S)=A^+$, the condition $\phi (U')\in A(S)$ for $U'\in B^+$ is vacuously always true. So Property (P6) takes the following form.

\begin{lemma}\label{S_condition_P6} $U'\to_{R_S}^*\rho(\phi (U'))$ for all $U'\in B^+$. \textnormal{(Property (P6)).}
\end{lemma}

\begin{proof} Let $U'\in B^+$. We shall prove by induction on $\Vert U'\Vert$ that $U'\to_{R_S}^* \rho(\phi(U'))$. Suppose $\Vert U'\Vert=1$. Then $U'=a$ for some $a\in A$ or $U'=s$  (recall that $B=A\cup \{s\}$). In either cases, we have $\rho(\phi(U'))=U'$. So $U'\to_{R_S}^* \rho(\phi(U'))$.

Suppose $\Vert U'\Vert>1$. Assume the lemma holds for all $U_1'\in B^+$ with  $\Vert U_1'\Vert<\Vert U'\Vert$. 

\noindent
{\bf Case 1.} Suppose $U'\in A^+$. Then $\phi (U')=U'$. If $U'$ ends with the subword $W_0$, say $U'=X_1W_0$ for some $X_1\in A^*$, then $\rho (U')=\rho (X_1)s=\rho (\phi (X_1))s$.  Since $W_0\to s\in R_S$ (the rule of the form ($\mathcal C$2)), we see that $U'\to_{R_S} X_1s$. Clearly $\Vert X_1\Vert <\Vert U'\Vert$. So by induction, $X_1\to_{R_S}^* \rho (\phi (X_1))$. Thus $U'=X_1W_0\to_{R_S} X_1s\to_{R_S}^* \rho (\phi (X_1))s=\rho (\phi(U'))$.

If $U'$ does not end with the subword $W_0$, then $U'=X_2a$ for some $X_2\in A^+$ and $a\in A$. Note that $\rho (U')=\rho (X_2)a=\rho (\phi(X_2))a$. By induction, $X_2\to_{R_S}^* \rho (\phi (X_2))$. Thus $U'\to_{R_S}^*\rho (\phi (U'))$.

\noindent
{\bf Case 2.} Suppose $U'=U_1'sU_2'$ for some $U_2'\in A^+$ and $U_1'\in B^*$. Note that $\phi (U')=\phi (U_1')W_0U_2'$. If $U_2'$ ends with the subword $W_0$, say $U_2'=X_1W_0$ for some $X_1\in A^*$, then $\rho (\phi (U'))=\rho (\phi (U_1')W_0X_1W_0)=\rho (\phi (U_1')W_0X_1)s=\rho (\phi (U_1'sX_1))s$. By induction, $U_1'sX_1\to_{R_S}^*\rho (\phi (U_1'sX_1))$. Also $U'\to_{R_S} U_1'sX_1s$ by the rule $W_0\to s\in R_S$ (rule ($\mathcal C$2)). 
Thus $U'\to_{R_S}^*\rho (\phi (U'))$.

Suppose $U_2'$ does not end with the subword $W_0$, but $W_0U_2'$ ends with the subword $W_0$, say $W_0U_2'=X_2W_0$ for some $X_2\in A^+$. Then there is a $X_3\in A^+$ such that $W_0=X_2X_3=X_3U_2'$. So  $sU_2'\to X_2s\in R_S$ (a rule of the form ($\mathcal C$6)) and $U'\to_{R_S} U_1'X_2s$. On the other hand, $\rho (\phi (U'))=\rho (\phi (U_1')X_2W_0)=\rho (\phi (U_1')X_2)s=\rho (\phi (U_1'X_2))s$, and also $\Vert U_1'X_2\Vert=\Vert U_1'\Vert+\Vert X_2\Vert=\Vert U_1'\Vert+\Vert U_2'\Vert<\Vert U'\Vert$. Therefore by induction, $U_1'X_2\to_{R_S}^*\rho (\phi (U_1'X_2))$.  Thus $U'\to_{R_S} U_1'X_2s\to_{R_S}^*\rho (\phi (U_1'X_2))s=\rho (\phi (U'))$.

Suppose $W_0U_2'$ does not end with the subword $W_0$. Let $U_2'=U_3'a$ for some $a\in A$ and $U_3'\in A^*$. Note that $\rho (\phi (U'))=\rho (\phi (U_1')W_0U_3'a)=\rho (\phi (U_1')W_0U_3')a=\rho (\phi (U_1'sU_3'))a$. By induction $U_1'sU_3'\to_{R_S}^* \rho (\phi (U_1'sU_3'))$. Thus $U'\to_{R_S}^* \rho (\phi (U'))$.

\noindent
{\bf Case 3.} Suppose $U'=U_1's$ for some $U_1'\in B^+$. Note that $\phi (U')=\phi (U_1')W_0$ and $\rho (\phi (U'))=\rho(\phi (U_1'))s$. By induction, $U_1'\to_{R_S}^*\rho (\phi (U_1'))$, and thus $U'\to_{R_S}^*\rho (\phi (U'))$.

The proof of this lemma is complete.
\end{proof}

Since $A(S)=A^+$, we have $\phi (U')\in A(S)$ for all $U'\in B^+$. Therefore the following lemma holds by choosing $U''=U'$.

\begin{lemma}\label{S_condition_P4} For each $U'\in B^+$ there is a $U''\in B^+$ such that $\phi(U'')\in A(S)$ and $U'\to_{R_S}^* U''$. \textnormal{(Property (P4)).}
\end{lemma}

\begin{lemma}\label{S_condition_pre_P3} Suppose $U'\to_{R_S} V'$ by one of the rules of the form ($\mathcal C$1), ($\mathcal C$3), ($\mathcal C$4) or ($\mathcal C$5). Then $\phi(U')\neq \phi (V')$.
\end{lemma}

\begin{proof} Note that all the rules ($\mathcal C$1), ($\mathcal C$3), ($\mathcal C$4) or ($\mathcal C$5)
have the form $\rho(X)\to \rho(Y)$ where $X\neq Y$ and $X\to_{R}^* Y$. 

Let $U'=Z_1'\rho(X)Z_2'$ with $Z_1',Z_2'\in B^*$. Then $V'=Z_1'\rho(Y)Z_2'$. Note that $\phi(U')=\phi(Z_1')X\phi(Z_2')$ and $\phi(V')=\phi(Z_1')Y\phi(Z_2')$ (by Lemma \ref{S_condition_P5} and the fact that $\phi$ is a homomorphism). If $\phi(U')=\phi(V')$, then $X=Y$ and 

\begin{equation}
X\to_R Y\to_R X\to_R Y\to_R\cdots,\notag
\end{equation}
would be an infinite reduction sequence, contrary to the fact that $R$ is complete. Hence $\phi(U')\neq \phi(V')$.
\end{proof}

\begin{lemma}\label{S_condition_P3} There does not exist an infinite reduction sequence 
\begin{equation}
U_1'\to_{R_S}U_2'\to_{R_S}U_3'\to_{R_S}\cdots,\notag
\end{equation}
 of words from $B^+$ such that $\phi(U_1')=\phi(U_2')=\phi(U_3')=\cdots$. \textnormal{(Property (P3)).}
\end{lemma}

\begin{proof} Suppose that such a sequence exists. Since $\phi(U_i')=\phi(U_{i+1}')$, by Lemma \ref{S_condition_pre_P3}, we conclude that $U'_i\to_{R_S} U'_{i+1}$ by one of the rules of the form ($\mathcal C$2) or  ($\mathcal C$6). Note that if a rule of the form ($\mathcal C$2) is applied  to $U'_i\to_{R_S} U'_{i+1}$, then $\Vert U'_{i+1}\Vert< \Vert U'_i\Vert$. If a rule of the form ($\mathcal C$6) is applied  to $U'_i\to_{R_S} U'_{i+1}$, then $\Vert U'_{i+1}\Vert=\Vert U'_i\Vert$ and one of the letter $s$ in $U'_{i+1}$ will be further to the right than it is in $U'_i$. Thus $\Vert U_i'\Vert\geq \Vert U_{i+1}'\Vert$ for all $i$.
 
There is an integer $i_0$ such that for all $i\geq i_0$, $\Vert U_i'\Vert=\Vert U_{i+1}'\Vert$. So the only rule that can be applied on $U'_i\to_{R_S} U'_{i+1}$ is a rule of the form ($\mathcal C$6). Since one of the letter $s$ in $U'_{i+1}$ will be further to the right than it is in $U'_i$, this process cannot go on indefinitely. We have obtained a contradiction. Hence the lemma holds.
\end{proof}

\begin{lemma}\label{S_condition_P2} For any $U',V'\in B^+$ with $U'\to_{R_S}^* V'$, we have $\phi(U')\to_{R}^* \phi(V')$. \textnormal{(Property (P2)).}
\end{lemma}

\begin{proof} It is sufficient to show $U'\to_{R_S} V'$ with $U',V'\in B^+$ implies that $\phi(U')\to_{R}^* \phi(V')$.

Suppose $U'\to_{R_S} V'$ by a rule of the form ($\mathcal C$1), say $\rho(X)\to \rho(Y)\in R_S$ where $X\to Y\in R$. Let $U'=Z_1'\rho(X)Z_2'$ with $Z_1',Z_2'\in B^*$. Then $V'=Z_1'\rho(Y)Z_2'$. By Lemma \ref{S_condition_P5}, $\phi(U')=\phi(Z_1')X\phi(Z_2')$ and $\phi(V')=\phi(Z_1')Y\phi(Z_2')$. Clearly $\phi(U')\to_R\phi(V')$ by the rule $X\to Y$.

Suppose $U'\to_{R_S} V'$ by a rule of the form ($\mathcal C$2). Let $U'=Z_1'W_0Z_2'$ with $Z_1',Z_2'\in B^*$. Then $V'=Z_1'sZ_2'$. By Lemma \ref{S_condition_P5}, $\phi(U')=\phi(Z_1')W_0\phi(Z_2')=\phi(V')$. Clearly $\phi(U')\to_R^*\phi(V')$.

Suppose $U'\to_{R_S} V'$ by a rule of the form ($\mathcal C$3), say $\rho(Z_1X_1X_2)\to \rho(Y')$, where $X_1X_2\to Y_1\in R$, $W_0=Z_1X_1$ and $Z_1X_1X_2\to_R^* Y'$ ($X_1,Y_1\in A^+$, $X_2,Z_1\in A^*$ and $Y'\in\textnormal{Irr}(R)$). Let $U'=Z_3'\rho(Z_1X_1X_2)Z_4'$ with $Z_3',Z_4'\in B^*$. Then $V'=Z_3'\rho(Y')Z_4'$. By Lemma \ref{S_condition_P5}, $\phi(U')=\phi(Z_3')Z_1X_1X_2\phi(Z_4')$ and $\phi(V')=\phi(Z_3')Y'\phi(Z_4')$. So $\phi(U')\to_R^*\phi(V')$, for $Z_1X_1X_2\to_R^* Y'$.

Similarly we can show that if  $U'\to_{R_S} V'$ by a rule of the form ($\mathcal C$4), ($\mathcal C$5) or ($\mathcal C$6), then $\phi(U')\to_R^*\phi(V')$. The proof of this lemma is complete.
\end{proof}

\begin{lemma}\label{S_condition_P1} For any $U\in A(S)$ and $V_1\in A^+$ with $U\to_R V_1$, there is a $U'\in B^+$ such that $U\to_R V_1\to_R^* \phi(U')$ and  $\rho(U)\to_{R_S} U'$. \textnormal{(Property (P1)).}
\end{lemma}

\begin{proof} Let $U\to_R V_1$ by a rule $X_2\to Y_2\in R$. Let $U=X_1X_2X_3$ where $X_1,X_3\in A^*$. Then $V_1=X_1Y_2X_3$

\noindent
{\bf Case 1.} Suppose $\rho(X_1X_2X_3)=\rho(X_1X_2)\rho(X_3)$.

\noindent
{\bf SubCase 1.1.}  Suppose $\rho(X_1X_2)=\rho(X_1)\rho(X_2)$. Then $\rho(U)=\rho(X_1)\rho(X_2)\rho(X_3)$ and also  $\rho(U)\to_{R_S} \rho(X_1)\rho(Y_2)\rho(X_3)$ by the rule $\rho (X_2)\to \rho(Y_2)\in R_S$ (a rule of the form ($\mathcal C$1)). Let $U'=\rho(X_1)\rho(Y_2)\rho(X_3)$.  By Lemma \ref{S_condition_P5},  $\phi(U')=X_1Y_2X_3=V_1$ and thus the lemma holds.

\noindent
{\bf SubCase 1.2.} Suppose  $\rho(X_1X_2)\neq \rho(X_1)\rho(X_2)$. By Lemma \ref{S_property_rho}, there are $Z_1,Z_4\in A^*$ and $Z_2,Z_3\in A^+$ with $X_1=Z_1Z_2$, $X_2=Z_3Z_4$ and $Z_2Z_3=W_0$ such that $\rho(X_1X_2)=\rho(Z_1)s\rho(Z_4)$. Note that $\rho (Z_2Z_3Z_4)\to \rho(Y')\in R_S$ where $Z_2Z_3Z_4\to_R^* Y'$ and $Y'\in\textnormal{Irr} (R)$ (a rule of the form ($\mathcal C$3)). Furthermore $\rho(Z_1Z_2Z_3Z_4)=\rho(X_1X_2)=\rho(Z_1)s\rho(Z_4)=\rho(Z_1Z_2Z_3)\rho(Z_4)$. So by Lemma \ref{S_property_rho_ending}, $\rho (Z_2Z_3Z_4)=\rho(Z_2Z_3)\rho(Z_4)=s\rho(Z_4)$. Therefore $\rho(X_1X_2)=\rho(Z_1)s\rho(Z_4)\to_{R_S} \rho(Z_1)\rho(Y')$ and 
\begin{equation}
\rho(U)=\rho(X_1X_2)\rho(X_3)\to_{R_S} \rho(Z_1)\rho(Y')\rho(X_3).\notag
\end{equation}

Let $U'=\rho(Z_1)\rho(Y')\rho(X_3)$. Then by Lemma \ref{S_condition_P5}, $\phi(U')=Z_1Y'X_3$. Note that $Z_2Z_3Z_4\to_R Z_2Y_2\to_R^* Y'$ (for $Y'\in\textnormal{Irr}(R)$). Therefore 
\begin{equation}
U=(Z_1Z_2)(Z_3Z_4)X_3\to_R V_1=(Z_1Z_2)Y_2X_3\to_R^* \phi(U'),\notag
\end{equation}
and thus the lemma holds.

\noindent
{\bf Case 2.} Suppose $\rho(X_1X_2X_3)\neq \rho(X_1X_2)\rho(X_3)$. By Lemma \ref{S_property_rho}, there are $Z_1,Z_4\in A^*$ and $Z_2,Z_3\in A^+$ with $X_1X_2=Z_1Z_2$, $X_3=Z_3Z_4$ and $Z_2Z_3=W_0$ such that $\rho(X_1X_2X_3)=\rho(Z_1)s\rho(Z_4)$. Since $W_0\in \textnormal {Irr} (R)$, we must have $\Vert Z_2\Vert <\Vert X_2\Vert$ (if not, then $X_2$ would be a subword of $W_0$ and $W_0\notin\textnormal{Irr}(R)$ because $X_2\to Y_2\in R$). Let $X_2=X_4Z_2$ for some $X_4\in A^+$. Then $Z_1=X_1X_4$.

\noindent
{\bf SubCase 2.1.}  Suppose $\rho(X_1X_4)=\rho(X_1)\rho(X_4)$. Note that $\rho (X_4Z_2Z_3)\to \rho(Y')\in R_S$ where $X_4Z_2Z_3\to_R^* Y'$ and $Y'\in\textnormal{Irr} (R)$ (a rule of the form ($\mathcal C$4)).  Furthermore $\rho(X_4Z_2Z_3)=\rho(X_4)s$ and $\rho(U)=\rho(X_1X_2X_3)=\rho(Z_1)s\rho(Z_4)=\rho(X_1X_4)s\rho(Z_4)=\rho(X_1)\rho(X_4)s\rho(Z_4)\to_{R_S} \rho(X_1)\rho(Y')\rho(Z_4)$. Let $U'=\rho(X_1)\rho(Y')\rho(Z_4)$. Then by Lemma \ref{S_condition_P5}, $\phi(U')=X_1Y'Z_4$. As before $X_4Z_2Z_3\to_R Y_2Z_3\to_R^* Y'$ (recall that $X_2=X_4Z_2$) and
\begin{equation}
U=(Z_1Z_2)(Z_3Z_4)=(X_1X_4Z_2)(Z_3Z_4)\to_R X_1Y_2Z_3Z_4=V_1\to_R^* \phi(U').\notag
\end{equation}
So the lemma holds.

\noindent
{\bf SubCase 2.2.}  Suppose $\rho(X_1X_4)\neq \rho(X_1)\rho(X_4)$. By Lemma \ref{S_property_rho}, there are $Z_5,Z_8\in A^*$ and $Z_6,Z_7\in A^+$ with $X_1=Z_5Z_6$, $X_4=Z_7Z_8$ and $Z_6Z_7=W_0$ such that $\rho(X_1X_4)=\rho(Z_5)s\rho(Z_8)$. Note that 
\begin{equation}
U=X_1X_2X_3=Z_5Z_6(Z_7Z_8Z_2)Z_3Z_4,\notag
\end{equation}
and $X_2=Z_7Z_8Z_2$. Also $\rho (Z_6(Z_7Z_8Z_2)Z_3)\to \rho(Y')\in R_S$ where $Z_6(Z_7Z_8Z_2)Z_3\to_R^* Y'$ and $Y'\in\textnormal{Irr} (R)$ (a rule of the form ($\mathcal C$5)).  Since $\rho(Z_5Z_6Z_7Z_8)=\rho(X_1X_4)=\rho(Z_5)s\rho(Z_8)=\rho(Z_5Z_6Z_7)\rho(Z_8)$, by Lemma \ref{S_property_rho_ending}, $\rho(Z_6Z_7Z_8)=\rho(Z_6Z_7)\rho(Z_8)=s\rho(Z_8)$. So $\rho (Z_6(Z_7Z_8Z_2)Z_3)=\rho(Z_6Z_7Z_8)s=s\rho (Z_8)s$ and $s\rho (Z_8)s\to \rho(Y')\in R_S$.

Recall that 
\begin{align}
\rho (Z_5Z_6(Z_7Z_8Z_2)Z_3Z_4) =\rho (U) & =\rho(X_1X_2X_3)\notag\\
&=\rho(Z_1)s\rho(Z_4)\notag\\
&=\rho(X_1X_4)s\rho(Z_4)\notag\\
&=\rho(Z_5)s\rho(Z_8)s\rho(Z_4).\notag
\end{align}

Therefore $\rho(U)=\rho(Z_5)s\rho(Z_8)s\rho(Z_4)\to_{R_S} \rho(Z_5)\rho(Y')\rho(Z_4)$. Let $U'=\rho(Z_5)\rho(Y')\rho(Z_4)$. Then by Lemma \ref{S_condition_P5}, $\phi(U')=Z_5Y'Z_4$. As before $Z_6(Z_7Z_8Z_2)Z_3\to_R Z_6Y_2Z_3\to_R^* Y'$ (recall that $X_2=X_4Z_2=Z_7Z_8Z_2$) and
\begin{equation}
U=Z_5Z_6(Z_7Z_8Z_2)Z_3Z_4\to_R Z_5Z_6Y_2Z_3Z_4=V_1\to_R^* \phi(U').\notag
\end{equation}
The proof of this lemma is complete.
\end{proof}

By Corollary \ref{general_propert_presentation_S},  Lemma \ref{S_condition_P1}, Lemma \ref{S_condition_P2}, Lemma \ref{S_condition_P3}, Lemma \ref{S_condition_P4}, Lemma \ref{S_condition_P5} and Lemma \ref{S_condition_P6}, we have shown that $[B\ ;\ R_S]$ is a semigroup presentation for $S$, $R_S$ is a finite complete rewriting system  and $W_0\to_{R_S}^* s$. Now note that if $U'\to V'\in R_S$ is a rule of the form ($\mathcal C$2), ($\mathcal C$3), ($\mathcal C$4), ($\mathcal C$5) or ($\mathcal C$6), then $\Vert U'\Vert>1$. From this we conclude that $s\in \textnormal {Irr} (R_S)$. Note also that if $X\in A^+$, $X\neq W_0$ and $\Vert X\Vert>1$, then $\Vert \rho(X)\Vert>1$. Therefore if $X\to Y\in R$ with $\Vert X\Vert>1$, then $\rho(X)\to \rho (Y)\in R_S$ and $\Vert \rho(X)\Vert>1$ (a rule of the form ($\mathcal C$1)). This implies that if $a\in A\cap \textnormal {Irr} (R)$, then $a\in \textnormal {Irr} (R_S)$.

Thus we have proved the following theorem.

\begin{theorem}\label{S_change_presentation} Let $[ A\ ;\ R]$ be a finitely presented semigroup presentation for $S$ for which $R$ is complete. Let $W_0\in A^+$ be such that $\Vert W_0\Vert>1$ and $W_0\in \textnormal {Irr} (R)$. Now let $s$ be a symbol that does not appear in $A$ and set $B=A\cup \{s\}$. Then there is complete rewriting system $R_S$ such that $[B\ ;\ R_S]$ is a finitely presented semigroup presentation for $S$ and $W_0\to_{R_S}^* s$. Furthermore $s\in\textnormal{Irr} (R_S)$, and  $a\in \textnormal {Irr} (R_S)$ for all $a\in A\cap \textnormal {Irr} (R)$.
\end{theorem}

\section {\bf Reduction process.}

In this section we will make further refinements and improvements (we call them reductions) to Theorem \ref{S_change_presentation}. The reason for such reductions is that we need a finitely presented semigroup presentation for $S$, which can be handled easily.

Let $S$ be a semigroup and $T$ be a large subsemigroup of $S$. Let $[ A\ ;\ R]$ be a finitely presented semigroup presentation for $S$ for which $R$ is complete. Let $S\setminus T=\{ [W_1]_R, [W_2]_R,\dots ,[W_n]_R\}$ with $W_i\in\textnormal {Irr} (R)$ and $\Vert W_1\Vert\leq \Vert W_2\Vert\leq\cdots \leq \Vert W_n\Vert$. Suppose that $\Vert W_1\Vert=\Vert W_2\Vert=\cdots =\Vert W_{i_0-1}\Vert=1$ and $\Vert W_{i_0}\Vert>1$. By Theorem \ref{S_change_presentation}, there is a finitely presented semigroup presentation $[B_{i_0}\ ;\ R_{i_0}]$ for $S$ such that $B=A\cup \{s_{i_0}\}$ for some symbol $s_{i_0}$  that does not appear in $A$, $R_{i_0}$ is complete, $W_{i_0}\to_{R_{i_0}}^* s_{i_0}$ and $W_1,W_2,\dots,  W_{i_0-1}, s_{i_0}\in\textnormal {Irr} (R_{i_0})$.

Now in this new semigroup presentation $[B_{i_0}\ ;\ R_{i_0}]$, we see that 

\begin{equation}
S\setminus T=\{ [W_1]_{R_{i_0}}, [W_2]_{R_{i_0}},\dots ,[W_{i_0-1}]_{R_{i_0}}, [s_{i_0}]_{R_{i_0}},[W_{i_0+1}']_{R_{i_0}},\dots ,[W_n']_{R_{i_0}}\},\notag 
\end{equation}
with $W_1,\dots, W_{i_0-1}, s_{i_0}, W_{i_0+1}',\dots ,W_{n}'\in\textnormal {Irr} (R_{i_0})$. 

Note that this process can be continued (in at most $n$ steps) until we obtain a finitely presented semigroup presentation $[B_{n}\ ;\ R_{n}]$ for $S$ such that $R_n$ is complete and $S\setminus T=\{ [s_1]_{R_{n}}, [s_2]_{R_{n}},\dots ,[s_n]_{R_{n}}\}$
with $s_1,\dots, s_{n}\in\textnormal {Irr} (R_{n})\cap B_n$. 

In fact by a standard procedure described in \cite[Section 2.2]{BO}, we may further assume that
for each $X\to Y\in R_n$, we have $Y\in\textnormal {Irr} (R)$, and for each $X\to Y\in R_n$, there is no $X'\in B_n^+$ for which $X\to_{R_n} X'$ by any rule in $R_n\setminus\{X\to Y\}$. This is the form of the presentation that we will use.

\section {\bf The main result.}

Let $S$ be a semigroup and $T$ be a large subsemigroup of $S$. As stated in Section 4, we may assume that $[ A\ ;\ R]$ is a finitely presented semigroup presentation for $S$ for which $R$ is  complete  and 
\begin{itemize}
\item [(Q1)] $S\setminus T=\{ [s_1]_{R}, [s_2]_{R},\dots ,[s_n]_{R}\}$ with $s_1,\dots, s_{n}\in\textnormal {Irr} (R)\cap A$,
\item [(Q2)] for each $X\to Y\in R$, we have $Y\in\textnormal {Irr} (R)$,
\item [(Q3)] for each $X\to Y\in R$, there is no $X'\in A^+$ for which $X\to_{R} X'$ by any rule in $R\setminus\{X\to Y\}$.
\end{itemize}

In order to show that $T$ has a finite complete rewriting system, we shall find a 5-tuple $(B,R_T,A(T),\phi, \rho)$ that has Property $\mathcal R$ relative to $[ A\ ;\ R]$ and apply Theorem \ref{general_propert_R}.

Let $A_1=\{ a\in A\ :\ [a]_R\in T\}$ and $A_S=\{s_1,s_2,\dots ,s_n\}$. Note that in general the union of $A_S$ and $A_1$ is not necessary equal to $A$. This is because there might exist an element $b\in A$ such that $[b]_R\in S\setminus T$. If this happens, we would have $b\to_R^* s_i$ for some $i$. 

\begin{lemma}\label{element_in_R} Let $X\to Y\in R$ with $[X]_R\in T$. Then 
\begin{itemize}
\item [(a)] if $W\in A^+$ is a subword of $X$ and $[W]_R\in S\setminus T$, then $W=s_i$ for some $i$,
\item [(b)] if $W\in A^+$ is a subword of $Y$ and $[W]_R\in S\setminus T$, then $W=s_i$ for some $i$.
\end{itemize}
\end{lemma}

\begin{proof} (a) Suppose $W\notin A_S$. Then by (Q1) $W\to_R^* s_i$ for some $i$. To be precise there is a $W_1\in A^+$ such that $W\to_R W_1\to_R^* s_i$. Let $W\to_R W_1$ by the rule $X_1\to Y_1$. Since $[X]_R\in T$, we cannot have $W=X$. Therefore $X_1\neq X$ and $X_1\to Y_1\in R\setminus \{X\to Y\}$. Let $X=Z_1WZ_2$ where $Z_1,Z_2\in A^*$. Then $X\to_R Z_1W_1Z_2$ by the rule $X_1\to Y_1$, contrary to (Q3). Hence $W=s_i$ for some $i$.

(b) can be proved similarly using the fact that $Y\in \textnormal {Irr} (R)$ (see (Q2)).
\end{proof}

We now begin to define the 5-tuple $(B,R_T,A(T),\phi, \rho)$. Let $A(T)(0)$ be the set of all $W\in A^+$, such that $[W]_R\in T$, and if  $X_1$ is a subword $W$ with $[X_1]_R\in S\setminus T$, then $\Vert X_1\Vert=1$ and $X_1\in A_S$. In other word, 
\begin{align}
A(T)(0)=\{W\in (A_1\cup A_s)^+ \ :&\ [W]_R\in T, \ \textnormal{and $W$ does not contain any subword}\notag\\
&\hskip 1cm\textnormal{$X_1$ with  $[X_1]_R\in S\setminus T$ and $\Vert X_1\Vert>1$}\}.\notag
\end{align}
 The following lemma is clear from the definition of $A(T)(0)$. 

\begin{lemma}\label{subword_in_AT0} Let $W\in A(T)(0)$ and $W'$ be a subword of $W$. If $[W']_R\in T$, then  $W'\in A(T)(0)$.
\end{lemma}

Next let
\begin{align}
F_1 & =A_1,\notag\\
F_2 & =\{ sb\ :\ s\in A_S, b\in A_1\cup A_S\ \textnormal{and}\ [sb]_R\in T\},\notag\\
F_3 & =\{ as\ :\  a\in A_1, s\in A_S\ \textnormal{and}\ [as]_R\in T\},\notag\\
F_4 & =\{ sbs'\ :\  s,s'\in A_S, b\in A_1\cup A_S\ \textnormal{and}\ [sb]_R, [bs']_R, [sbs']_R\in T\}.\notag
\end{align}
It is not hard to see that if $W\in F_1\cup F_2\cup F_3\cup F_4$, then $[W]_R\in T$. Furthermore $F_1\cup F_2\cup F_3\cup F_4\subseteq A(T)(0)$. For convenience, for each $G\subseteq A^+$ and $X\in A^+$, we set $XG=\{ XW\ :\ W\in G\}$. 

Now we shall define $A(T)$. Let $A(T)(1)= F_1\cup F_2\cup F_3\cup F_4$ and for each $i\geq 1$, let
\begin{equation}
A(T)(i+1)=\left (\bigcup_{a\in A_1} \left( (aA(T)(i))\cap A(T)(0)\right)\right ) \cup \left (\bigcup_{X\in F_2} \left ((XA(T)(i))\cap A(T)(0)\right )\right ).\notag
\end{equation}
Set $A(T)=\bigcup_{i\geq 1} A(T)(i)$. In the following lemma we shall prove some properties of $A(T)$.
\begin{lemma}\label{property_A(T)} \
\begin{itemize}
\item [(a)] $A(T)=A(T)(0)$.
\item [(b)] $A(T)$ contains the set $\{ W\in \textnormal{Irr} (R)\ :\ [W]_R\in T\}$.
\item [(c)] Let $X\to Y\in R$ with $[X]_R\in T$. Then $X,Y\in A(T)$. 
\end{itemize}
\end{lemma}

\begin{proof} (a) Clearly $A(T)\subseteq A(T)(0)$. Let $W\in A(T)(0)$. We shall prove by induction on $\Vert W\Vert$ that $W\in A(T)$. 

Suppose $\Vert W\Vert=1$. Since $[W]_R\in T$, we must have $W\in A_1$. So $W\in A(T)(1)\subseteq A(T)$. 

Suppose $\Vert W\Vert=2$. Then $W=a'a$, or $W=as$, or $W=sa$, or $W=ss'$ ($a,a'\in A_1$, $s,s'\in A_S$). If $W=a'a$, then $W\in (a'A(T)(1))\cap A(T)(0)\subseteq A(T)(2)\subseteq A(T)$. If $W=as$, then $W\in F_3\subseteq A(T)(1)\subseteq A(T)$. If $W=sa$ or $W=ss'$, then $W\in F_2\subseteq A(T)(1)\subseteq A(T)$. 

Suppose $\Vert W\Vert \geq 3$. Assume that it is true for all $W'\in A(T)(0)$ with $\Vert W'\Vert<\Vert W\Vert$.

If $W$ begins with a letter $a\in A_1$, say $W=aW'$ where $W'\in A^+$, then $\Vert W'\Vert\geq 2$. Note that $[W']_R\in T$, for if $[W']_R\in S\setminus T$, then by the definition of $A(T)(0)$, $W'\in A_S$ and $\Vert W'\Vert=1$, contrary to the fact that $\Vert W'\Vert\geq 2$. Therefore by Lemma \ref{subword_in_AT0},  $W'\in A(T)(0)$. By induction, $W'\in A(T)$. Let $W'\in A(T)(i)$ for some $i\geq 1$. Then $W\in (aA(T)(i))\cap A(T)(0)\subseteq A(T)(i+1)\subseteq A(T)$.

If $W$ begins with a letter $s\in A_S$, say $W=sbW'$ where $b\in A_1\cup A_S$ and $W'\in A^+$, then $\Vert W'\Vert\geq 1$. If $[W']_R\in S\setminus T$, then by the definition of $A(T)(0)$, $W'=s'$ for some $s'\in A_S$, and $W=sbs'$. Since $W\in A(T)(0)$, we have $[sb]_R, [bs']_R, [sbs']_R\in T$  (definition of $A(T)(0)$). This means $W\in F_4\subseteq A(T)(1)\subseteq A(T)$.

If $[W']_R\in T$, then by Lemma \ref{subword_in_AT0},  $W'\in A(T)(0)$. By induction, $W'\in A(T)$. Let $W'\in A(T)(i)$ for some $i\geq 1$. Then $W\in (sbA(T)(i))\cap A(T)(0)\subseteq A(T)(i+1)\subseteq A(T)$.

The proof of part (a) of the lemma is complete.

Part (b) follows from part (a) and the fact that $A(T)(0)$ contains the set $\{ W\in \textnormal{Irr} (R)\ :\ [W]_R\in T\}$.

\noindent
(c) By part (a) of Lemma \ref{element_in_R}, we conclude that $X$ does not contain any subword $X_1$ with $[X_1]_R\in S\setminus T$ and $X_1\notin A_S$. So $X\in A(T)(0)=A(T)$. Similarly by part (b) of Lemma \ref{element_in_R},  $Y\in A(T)$.
\end{proof}

Now we shall define the set $B$ and the homomorphism $\phi$. Let

\begin{align}
C_R &=\{ c_{as}\ :\ [as]_R\in T\ \textnormal{with $a\in A_1$ and $s\in A_S$}\},\notag\\
C_{L_1} &=\{ c_{sa}\ :\ [sa]_R\in T\ \textnormal{with $a\in A_1$ and $s\in A_S$}\},\notag\\
C_{L_2} &=\{ c_{ss'}\ :\ [ss']_R\in T\ \textnormal{with $s,s'\in A_S$}\},\notag\\
C_{M_1} &=\{ c_{s'as}\ :\ [s'as]_R, [s'a]_R,[as]_R\in T\ \textnormal{with $a\in A_1$ and $s,s'\in A_S$}\},\notag\\
C_{M_2} &=\{ c_{ss's''}\ :\ [ss's'']_R, [ss']_R,[s's'']_R\in T\ \textnormal{with $s,s',s''\in A_S$}\}.\notag
\end{align}

Set $C=C_R\cup C_{L_1}\cup C_{L_2}\cup C_{M_1}\cup C_{M_2}$ and  $B=A_1\cup C$. Since $A_1$ and $A_S$ are finite, it is not hard to see that $B$ is finite. Let $\phi_1: B\to A^+$ be defined by $\phi_1(a)=a$ for all $a\in A_1$ and $\phi_1(c_{u})=u$ for all $c_{u}\in C$ (for example $\phi_1(c_{as})=as$ for $c_{as}\in C_R$). Clearly $\phi_1$ can be extended to a homomorphism $\phi : B^+\to A^+$ by defining $\phi (U')=\phi_1(b_1)\dots \phi_1(b_l)$ for all $U'=b_1\dots b_l\in B^+$. Furthermore $[\phi(U')]_R\in T$ for all $U'\in B^+$. For convenience, we may define $\phi(\epsilon_B)=\epsilon_A$ where $\epsilon_B$ and $\epsilon_A$ are empty words in $B^*$ and $A^*$, respectively. 
The following lemma is obvious.

\begin{lemma}\label{T_property_phi} For all $U'\in B^+$, $\Vert \phi(U')\Vert\geq \Vert U'\Vert$.
\end{lemma}

We define $\rho : A(T)\to B^+$ as follows:

Let $W\in A(T)$. 
\begin{itemize}
\item [(a)] Suppose $W\in A(T)(1)$. If $W\in F_1$, then set $\rho (W)=W$. If $W\in F_2\cup F_3\cup F_4$, set $\rho (W)=c_W$ (for example if $W=as\in F_3$, then $\rho(W)=c_{as}$).
\item [(b)] Suppose $W\in A(T)(i+1)$ for some $i\geq 1$. Then $W=aW_1$ or $W=sbW_1$ ($a\in A_1$, $s\in A_S$, $b\in A_1\cup A_S$ and $W_1\in A(T)(i)$). If the former holds, set $\rho(W)=a\rho(W_1)$. If the latter holds, set $\rho(W)=c_{sb}\rho(W_1)$.
\end{itemize}
The function $\rho$ is well-defined can be easily proved by observing that a word from $A(T)(i+1)$ is obtained in a unique way from a unique word from $A(T)(i)$. As for the homomorphism $\phi$,  we may define $\rho(\epsilon_A)=\epsilon_B$.

\begin{lemma}\label{T_property_rho_ending} Let $U\in A(T)(l)$ for some $l\geq 1$. Then $\rho(U)=b_1'\dots b_l'$ where $b_i'\in B$. Furthermore if $l>1$, then $b_i'\in A_1\cup C_{L_1}\cup C_{L_2}$ for all $1\leq i\leq l-1$.
\end{lemma}
\begin{proof} We prove by induction on $l$. Suppose $l=1$. Then $\rho(U)=b_1'$ by the definition of $\rho$. Suppose $l>1$. Assume that it is true for all $l'$ with $l'<l$.

Since $U\in A(T)(l)$, we have either $U=aU_1$ or $U=sbU_1$ ($a\in A_1$, $sb\in F_2$ and $U_1\in A(T)(l-1)$).
Suppose $U=aU_1$. Then $\rho(U)=a\rho(U_1)$. This means $b_1'=a\in A_1$. By induction $\rho(U_1)=b_2'\dots b_l'$. Furthermore if $l-1>1$ (i.e. $l>2$), then $b_2',\dots, b_{l-1}'\in A_1\cup C_{L_1}\cup C_{L_2}$.

Suppose $U=sbU_1$. Then $\rho(U)=c_{sb}\rho(U_1)$. This means $b_1'=c_{sb}\in C_{L_1}\cup C_{L_2}$. By induction $\rho(U_1)=b_2'\dots b_l'$. Furthermore if $l-1>1$ (i.e. $l>2$), then $b_2',\dots, b_{l-1}'\in A_1\cup C_{L_1}\cup C_{L_2}$.

Hence in either cases the lemma holds.
\end{proof}

\begin{lemma}\label{T_condition_P5} $\phi(\rho(U))=U$ for all $U\in A(T)$. \textnormal{(Property (P5)).}
\end{lemma}
\begin{proof} We just need to show that for all $i\geq 1$, if $U\in A(T)(i)$, then $\phi(\rho(U))=U$. 

Suppose $U\in A(T)(1)$. If $U\in  F_1$, then $\rho(U)=U$ and $\phi(\rho(U))=U$. If $U\in F_2\cup F_3\cup F_4$, then $\rho(U)=c_U$ and $\phi(\rho(U))=\phi(c_U)=U$. Assume that it is true for all $U'\in A(T)(i)$.

Let $U\in A(T)(i+1)$. Then $U=aU_1$ or $U=sbU_1$ where $a\in A_1$, $sb\in F_2$ and $U_1\in A(T)(i)$. If the former holds, then $\rho(U)=a\rho(U_1)$ and by induction $\phi(\rho(U))=a\phi(\rho(U_1))=aU_1=U$. If the latter holds, 
then $\rho(U)=c_{sb}\rho(U_1)$, and by induction $\phi(\rho(U))=\phi(c_{sb})\phi(\rho(U_1))=sbU_1=U$. Hence the lemma holds.
\end{proof}

\begin{lemma}\label{T_property_rho_equal} Let $U'=b_1'\dots b_l'\in B^+$ where $b_i'\in A_1\cup C_{L_1}\cup C_{L_2}$ for all $1\leq i\leq l-1$ and $b_l'\in B$. If $\phi(U')\in A(T)$, then $\phi(U')\in A(T)(l)$ and $\rho(\phi(U'))=U'$.
\end{lemma}
\begin{proof} We prove by induction on $l$. Suppose $l=1$. If $b_1'=a\in A_1$, then $\phi(b_1')=a$, and $\rho(\phi(b_1'))=b_1'$. If $b_1'=c_z\in C$, then $\phi(b_1')=z\in A(T)(1)$, and $\rho(\phi(b_1'))=b_1'$.

Suppose $l>1$. Assume that it is true for all $l'$ with $l'<l$. Let $U'=b_1'U_1'$ where $U_1'=b_2'\dots b_l'$.
By induction, $\phi(U_1')\in A(T)(l-1)$ and $\rho(\phi(U'_1))=U'_1$. Since $b_1'\in  A_1\cup C_{L_1}\cup C_{L_2}$, we have $\phi(b_1')\in A_1\cup F_2$. Therefore $\phi(U')=\phi(b_1')\phi(U_1')\in A(T)(l)$, and $\rho(\phi(U'))=b_1'\rho(\phi(U'_1))=b_1'U_1'=U'$. Hence the lemma holds.
\end{proof}

We are now ready to define the rules in $R_T$. Let us begin by recalling some of the results  of Lemma \ref{property_A(T)}. For each $X\to Y\in R$ with $[X]_R\in T$, we have $X,Y\in A(T)$ (part (c) of Lemma \ref{property_A(T)}). Furthermore if $Y\in \textnormal{Irr} (R)$ and $[Y]_R\in T$, then $Y\in A(T)$ (part (b) of Lemma \ref{property_A(T)}). Recall that $C=C_R\cup C_{L_1}\cup C_{L_2}\cup C_{M_1}\cup C_{M_2}$, $\epsilon_A$ is the empty word in $A^*$, $\phi$ is a homomorphism of $B^+$ into $A^+$ (furthermore $[\phi(U')]_R\in T$ for all $U'\in B^+$), and $\rho$ is a function of  $A(T)$ into $B^+$. As $R$ is a finite complete rewriting system, $\textnormal {Left} (R)=\{ X\in A^+\ :\ X\to Y\in R\}$ is finite. Let $N=\left(\max_{X\in \textnormal {Left} (R)} \Vert X\Vert\right)+4$. 
The rules are grouped into two forms, ($\mathcal D$1) and ($\mathcal D$2):

\begin{itemize}
\item [($\mathcal D$1)] for each $U'\in B^+$ with $\Vert \phi(U')\Vert\leq N$ and $\phi(U')\notin \textnormal{Irr} (R)$, put $U'\to \rho(\overline {\phi(U')})$ in $R_T$ where $\phi(U')\to_R^*\overline {\phi(U')}$ and $\overline {\phi(U')}\in \textnormal{Irr} (R)$;
\item [($\mathcal D$2)] for each $U'\in B^+$ with $\Vert U'\Vert=2$, $\phi(U')\in A(T)$ and $U'\neq \rho(\phi(U'))$, put $U'\to \rho(\phi(U'))$ in $R_T$.
\end{itemize}

Note that the number of rules of the form ($\mathcal D$1) that we put in $R_T$ is finite, for by Lemma \ref{T_property_phi} the length of $U'$ is bounded and $B$ is finite. Similarly the number of rules of the form ($\mathcal D$2) that we put in $R_T$ is also finite. Therefore $R_T$ is  finite and $[B\ ;\ R_T]$ is finitely presented.  Note that by the main result in \cite[Theorem 6.1]{Ruskuc}, one can get a finite presentation for $T$ by taking $N$ sufficiently large.

\begin{lemma}\label{T_condition_Pre2_P6} Let $U',V'\in B^+$. If $U'\to_{R_T} V'$ by a rule of the form ($\mathcal D$2), then $\phi(U')=\phi(V')$. Furthermore either 
\begin{itemize}
\item[(i)] the number of elements in $C_R\cup C_{M_1}\cup C_{M_2}$ which appear as letters in the word $V'$ is less than that in the word $U'$, or 
\item[(ii)] the number of elements in $C_R\cup C_{M_1}\cup C_{M_2}$ which appear as letters in the word $V'$ is the same as that in the word $U'$, $\Vert U'\Vert=\Vert V'\Vert$, and  there is an element in $C_R\cup C_{M_1}\cup C_{M_2}$ in which it ``moves'' further right in the resulting word $V'$ than it is in the word $U'$ (the element may have changed).
\end{itemize}
\end{lemma}

\begin{proof} Let $U'\to_{R_T} V'$ by the rule $X'\to \rho(\phi(X'))$ where $X'\in B^+$, $\Vert X'\Vert =2$, $\phi(X')\in A(T)$ and $X'\neq \rho(\phi(X'))$. By Lemma \ref{T_condition_P5}, $\phi(\rho(\phi(X')))=\phi(X')$. Since $\phi$ is a homomorphism, we have $\phi(U')=\phi(V')$. Now we will show that either (i) or (ii) holds. 

If the first letter that appears in $X'$ is not from $C_R\cup C_{M_1}\cup C_{M_2}$, then by Lemma \ref{T_property_rho_equal}, $\rho(\phi(X'))=X'$, a contradiction. So we may assume that the first letter that appears in $X'$ is  from $C_R\cup C_{M_1}\cup C_{M_2}$.

By Lemma \ref{T_property_rho_ending}, $\rho(\phi(X'))$ has at most one letter from $C_R\cup C_{M_1}\cup C_{M_2}$, which is then the last letter.  If $\rho(\phi(X'))$ has no letter from $C_R\cup C_{M_1}\cup C_{M_2}$, then (i) holds. 

Suppose  $\rho(\phi(X'))$ has a letter from $C_R\cup C_{M_1}\cup C_{M_2}$. Then  $\phi(X')=\phi(\rho(\phi(X')))$ ends with a letter from $A_S$. Let $X'=cy$ where $c\in C_R\cup C_{M_1}\cup C_{M_2}$ and $y\in B$.  
Then $y\notin A_1\cup C_{L_1}$. If $y\in C_R\cup C_{M_1}\cup C_{M_2}$, then (i) holds. So we may assume that $y\in C_{L_2}$. Let $y=c_{s'''s''''}$. If $c=c_{as}$, then $\rho(\phi(X'))=ac_{ss'''s''''}$, if $c=c_{sas'}$, then $\rho(\phi(X'))=c_{sa}c_{s's'''s''''}$, and if $c=c_{ss's''}$, then $\rho(\phi(X'))=c_{ss'}c_{s''s'''s''''}$. 
Therefore $\Vert \rho(\phi(X'))\Vert=\Vert X'\Vert$ and (ii) holds.
\end{proof}

\begin{lemma}\label{T_condition_P6} $U'\to_{R_T}^*\rho(\phi (U'))$ for all $U'\in B^+$ with $\phi (U')\in A(T)$. \textnormal{(Property (P6)).}
\end{lemma}

\begin{proof} Let $U'=b_1'\dots b_l'\in B^+$ where $b_i'\in B$ for all $i$. If $b_i'\in A_1\cup C_{L_1}\cup C_{L_2}$ for all $1\leq i\leq l-1$, then by Lemma \ref{T_property_rho_equal}, $\rho(\phi(U'))=U'$. Hence $U'\to_{R_T}^*\rho(\phi (U'))$. 

So we may assume that $b_i'\in C_R\cup C_{M_1}\cup C_{M_2}$ for some $1\leq i\leq l-1$. By Lemma \ref{subword_in_AT0} and part (a) of Lemma \ref{property_A(T)}, $\phi(b_i'b_{i+1}')\in A(T)$. By Lemma \ref{T_property_rho_ending}, $b_i'b_{i+1}'\neq\rho(\phi(b_i'b_{i+1}'))$. Therefore $b_i'b_{i+1}'\to \rho(\phi(b_i'b_{i+1}'))$ is a rule of the form ($\mathcal D$2) in $R_T$. 

Let $V'=b_1'\dots b_{i-1}'\rho(\phi(b_i'b_{i+1}'))b_{i+2}'\dots b_l'$. Then $U'\to_{R_T} V'$, and by  Lemma \ref{T_condition_P5}, $\phi(U')=\phi(b_1'\dots b_l')=\phi(V')$. By  Lemma \ref{T_condition_Pre2_P6}, we conclude that after applying rules of the form ($\mathcal D$2) a finite number of times, there is a $U''=d_1'\dots d_r'\in B^+$ with $d_i'\in A_1\cup C_{L_1}\cup C_{L_2}$ for all $1\leq i\leq r-1$ and $d_r'\in B$, such that $U'\to_{R_T}^* U''$ and $\phi(U')=\phi(U'')$. Again by Lemma \ref{T_property_rho_equal}, $\rho(\phi(U''))=U''$. So $U'\to_{R_T}^* \rho(\phi(U''))=\rho(\phi(U'))$.
\end{proof}

\begin{lemma}\label{T_condition_Pre_P4} Let $U'\in B^+$ and $V\in A^+$. If $\phi(U')\to_R V$,  
 then there is a $V'\in B^+$ such that $U'\to_{R_T} V'$ by a rule of the form  \textnormal{(}$\mathcal D$1\textnormal{)}, and  $V\to_R^* \phi(V')$.
\end{lemma}

\begin{proof} Let $U'=b_1'\dots b_l'$ where $b_i'\in B$, and $\phi(U')\to_R V$ by a rule $X\to Y$ in $R$. Then for some non-negative integers $j_1,j_2$, $X$ is a subword of $\phi(b_{j_1}'\dots b_{j_1+j_2}')$. We may assume that $X$ is not a subword of $\phi(b_{j_1+1}'\dots b_{j_1+j_2}')$ or $\phi(b_{j_1}'\dots b_{j_1+j_2-1}')$. Since $\phi(b_{j_1}')$ and $\phi(b_{j_1+j_2}')$ are at most of length 3, we deduce that $\Vert \phi(b_{j_1}'\dots b_{j_1+j_2}')\Vert\leq \Vert X\Vert+4\leq N$. So $b_{j_1}'\dots b_{j_1+j_2}'\to \rho(\overline{\phi(b_{j_1}'\dots b_{j_1+j_2}')})$ is a rule of the form ($\mathcal D$1) in $R_T$, where $\phi(b_{j_1}'\dots b_{j_1+j_2}')\to_R^* \overline{\phi(b_{j_1}'\dots b_{j_1+j_2}')}$, $\phi(b_{j_1}'\dots b_{j_1+j_2}')\notin \textnormal{Irr} (R)$ and $\overline{\phi(b_{j_1}'\dots b_{j_1+j_2}')}\in \textnormal{Irr} (R)$. Set $V'=b_1'\dots b_{j_1-1}'\rho(\overline{\phi(b_{j_1}'\dots b_{j_1+j_2}')})b_{j_1+j_2+1}'\dots b_l'$. Then $U'\to_{R_T} V'$.

By Lemma \ref{T_condition_P5}, $\phi(V')=\phi(b_1'\dots b_{j_1-1}')(\overline{\phi(b_{j_1}'\dots b_{j_1+j_2}')})\phi(b_{j_1+j_2+1}'\dots b_l')$. Let $\phi(b_{j_1}'\dots b_{j_1+j_2}')=W_1XW_2$ where $W_1,W_2\in A^*$ (we allow $W_1, W_2$ to be empty word). Then $\phi(b_{j_1}'\dots b_{j_1+j_2}')\to_R W_1YW_2\to_R^*\overline{\phi(b_{j_1}'\dots b_{j_1+j_2}')}$. Hence $V=\phi(b_1'\dots b_{j_1-1}')(W_1YW_2)\phi(b_{j_1+j_2+1}'\dots b_l')\to_R^* \phi(V')$.
\end{proof}

\begin{lemma}\label{T_condition_P4} For each $U'\in B^+$ there is a $U''\in B^+$ such that $\phi(U'')\in A(T)$ and $U'\to_{R_T}^* U''$. \textnormal{(Property (P4)).}
\end{lemma}

\begin{proof} We shall prove by induction on $d_R(\phi(U'))$. Suppose $d_R(\phi(U'))=0$. Then $\phi(U')\in A(T)$  (part (b) of Lemma \ref{property_A(T)}). So we may choose $U''=U'$. Suppose $d_R(\phi(U'))>0$. Assume that it is true for all $U_1'\in B^+$ with $d_R(\phi(U_1'))<d_R(\phi(U'))$. 

Since $d_R(\phi(U'))>0$, there is a $V\in A^+$ such that $\phi(U')\to_R V$. By Lemma \ref{T_condition_Pre_P4}, there is a $V'\in B^+$ such that $U'\to_{R_T} V'$ and $V\to_R^* \phi(V')$.
Therefore $\phi(U')\to_R^* \phi(V')$, and $d_R(\phi(V'))<d_R(\phi(U'))$. By induction, there is a $U''\in B^+$ such that $\phi(U'')\in A(T)$ and $V'\to_{R_T}^* U''$. Hence $U'\to_{R_T}^* U''$.
\end{proof}

\begin{lemma}\label{T_condition_pre_P3} Suppose $U'\to_{R_T} V'$ by one of the rules of the form ($\mathcal D$1). Then $\phi(U')\neq \phi (V')$ and $\phi(U')\to_R^* \phi(V')$.
\end{lemma}

\begin{proof}  Suppose $U'\to_{R_T} V'$ by a rule of the form ($\mathcal D$1), say $X'\to Y'$. Then $\Vert\phi(X')\Vert\leq N$, $\phi(X')\notin \textnormal{Irr} (R)$, and $Y'=\rho(\overline {\phi(X')})$, where  $\phi(X')\to_R^*\overline {\phi(X')}$ and $\overline {\phi(X')}\in \textnormal{Irr} (R)$.

Let $U'=W_1'X'W_2'$ where $W_1', W_2'\in B^*$ (we allow $W_1'$ and $W_2'$ to be empty word). Note that $V'=W_1'\rho(\overline {\phi(X')})W_2'$. By Lemma \ref{T_condition_P5} and the fact that $\phi$ is a homomorphism, we must have $\phi(V')=\phi(W_1')\overline {\phi(X')}\phi(W_2')\neq \phi(U')$, for otherwise we would have $\phi(X')=\overline{\phi(X')}$. Furthermore $\phi(U')\to_R^* \phi(V')$.
\end{proof}

\begin{lemma}\label{T_condition_P3} There does not exist an infinite reduction sequence 
\begin{equation}
U_1'\to_{R_T}U_2'\to_{R_T}U_3'\to_{R_T}\cdots,\notag
\end{equation}
 of words from $B^+$ such that $\phi(U_1')=\phi(U_2')=\phi(U_3')=\cdots$. \textnormal{(Property (P3)).}
\end{lemma}

\begin{proof} Suppose that such a sequence exists. 

 Since $\phi(U_i')=\phi(U_{i+1}')$, by Lemma \ref{T_condition_pre_P3}, we conclude that $U_{i}'\to_{R_T}U_{i+1}'$ by a rule of the form ($\mathcal D$2). By Lemma \ref{T_condition_Pre2_P6}, the number of elements in $C_R\cup C_{M_1}\cup C_{M_2}$ which appear as letters in the word $U_{i+1}'$ is either less than that in the word $U_{i}'$, or the number are the same and $\Vert U_{i}'\Vert=\Vert U_{i+1}'\Vert$, but it `moves' to the right. So we deduce that there is an integer $i_0$ such that for all $i\geq i_0$, the number of elements in $C_R\cup C_{M_1}\cup C_{M_2}$ which appear as letters in the word $U_{i}'$ is the same as in the word  $U_{i+1}'$, and $\Vert U_{i}'\Vert=\Vert U_{i+1}'\Vert$.  So a letter (an element in $C_R\cup C_{M_1}\cup C_{M_2}$) in the word $U_{i}'$ will `move' further right in the word  $U_{i+1}'$. But this process cannot be continued indefinitely as $\Vert U_{i}'\Vert=\Vert U_{i+1}'\Vert$. We have obtained a contradiction.
\end{proof}

\begin{lemma}\label{T_condition_P2} For any $U',V'\in B^+$ with $U'\to_{R_T}^* V'$, we have $\phi(U')\to_{R}^* \phi(V')$. \textnormal{(Property (P2)).}
\end{lemma}

\begin{proof} It is sufficient to show $U'\to_{R_T} V'$ with $U',V'\in B^+$ implies that $\phi(U')\to_{R}^* \phi(V')$.

Suppose $U'\to_{R_T} V'$ by a rule of the form ($\mathcal D$1). By Lemma \ref{T_condition_pre_P3}, $\phi(U')\to_{R}^* \phi(V')$. Suppose $U'\to_{R_T} V'$ by a rule of the form ($\mathcal D$2). By Lemma \ref{T_condition_Pre2_P6}, $\phi(U')=\phi(V')$, and thus $\phi(U')\to_{R}^* \phi(V')$.
\end{proof}

\begin{lemma}\label{T_condition_P1} For any $U\in A(T)$ and $V_1\in A^+$ with $U\to_R V_1$, there is a $U'\in B^+$ such that $U\to_R V_1\to_R^* \phi(U')$ and  $\rho(U)\to_{R_T} U'$. \textnormal{(Property (P1)).}
\end{lemma}

\begin{proof} By Lemma \ref{T_condition_P5}, $U=\phi(\rho(U))$. By Lemma \ref{T_condition_Pre_P4}, there is a $U'\in B^+$ such that $\rho(U)\to_{R_T} U'$ by a rule of the form  \textnormal{(}$\mathcal D$1\textnormal{)}, and  $V_1\to_R^* \phi(U')$. The lemma follows. 
\end{proof}

\begin{proof}[\bf Proof of Theorem \ref{main_theorem}] Let $[ A\ ;\ R]$ be a finitely presented semigroup presentation for $S$ for which $R$ is complete. By the reduction process described in Section 4, we may assume that (Q1), (Q2) and (Q3) hold. Now the 5-tuple $(B,R_T,A(T),\phi, \rho)$ has been defined. By Theorem \ref{general_propert_R}, it is sufficient to show that $(B,R_T,A(T),\phi, \rho)$ has Property $\mathcal R$ relative to $[ A\ ;\ R]$. This has been done in Lemma \ref{T_condition_P5}, Lemma \ref{T_condition_P6}, Lemma \ref{T_condition_P4}, Lemma \ref{T_condition_P3}, Lemma \ref{T_condition_P2} and Lemma \ref{T_condition_P1}.
\end{proof}

\noindent
{\bf Acknowledgement.} We would like to thank Prof. S.J. Pride for introducing this problem to us. We are also indebted to the anonymous referee for suggesting condition (P3) and some of the proofs in Section 5, which helped us  improved this paper tremendously.

\end {document}